\documentclass[12pt]{amsart}
\numberwithin{equation}{section}
\newtheorem{ste}{Theorem}[section]
\newtheorem{lem}[ste]{Lemma}
\newtheorem{cor}[ste]{Corollary}
\newtheorem{prop}[ste]{Proposition}
\theoremstyle{definition}
\newtheorem{defi}[ste]{Definition}

\theoremstyle{remark}
\newtheorem{rem}[ste]{Remark}
 
\newenvironment{prettylist}{
\begin{list}{$\bullet$}{%
\settowidth{\labelwidth}{$\bullet$}%
\setlength{\leftmargin}{\labelsep}%
\addtolength{\leftmargin}{\labelwidth}%
}}{ \end{list}}
\usepackage{amssymb}
\usepackage{times}
\usepackage{graphicx}
\newcommand{\R}{{\mathbb R}}
\newcommand{\Co}[1]{{\mathcal C}^{#1}}
\DeclareMathOperator{\Emb}{Emb}
\DeclareMathOperator{\Image}{Image}
\DeclareMathOperator{\codim}{codim}
\DeclareMathOperator{\diffd}{d} % matica  DifferentialD
\begin{document}
\title{Applications of canonical relations in generic differential geometry}
\author{M. van Manen}
\address{Mathematisch Instituut \\ Universiteit Utrecht \\
P.\ O.\ box 80010 \\ 3508 TA Utrecht \\ The Netherlands}
\email{manen@math.uu.nl}
\date{September 12, 2002}
\begin{abstract}
Conflict sets are loci of intersecting wavefronts emanating from
$l$ different surfaces. We show that generically conflict sets
are Legendrian: locally they admit the structure of wavefronts.
Simple stable singularities for this problem in $\R^n$ occur
when $0 \leq n-l \leq 4$. Other related sets, such as kite curves
and centre sets are also defined 
and discussed. Throughout canonical relations are used as an
essential tool to carry out several of these geometrical
constructions.
\end{abstract}
\maketitle
\section*{Introduction}
The symmetry set of a manifold $M\subset \R^n$ 
is defined as
the closure of the set in $u \in \R^n$ where the distance function
\begin{equation*}
   s \overset{F}{\mapsto} \rVert u  -\gamma(s) \lVert 
\end{equation*}
has a double extremum.
Here $\gamma\colon M \rightarrow\R^n$ is an embedding of the
surface. The symmetry set is the closure
of the $A_1A_1$ stratum in the parameter space $\R^n$ of
the family of functions $F(x,s)$. 
\begin{figure}[htbp]
\centering
\includegraphics[width=8cm]{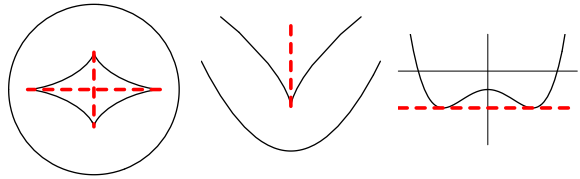}
\caption{Two symmetry sets and a double extremum}
\label{fig:1}
\end{figure}
The symmetry set was studied in \cite{BruceGiblinGibson}. There also 
a number of other related sets that measure symmetry  are discussed: medial
axis, cut-locus. For references to the many variants that 
exist we refer to this paper and other work of Giblin.
\vspace{1em}
\newline
In this paper we do not take one manifold, but $l$ manifolds of
codimension $1$ in $\R^n$. Hence, we have $l$ distance
functions. We are interested in the points where all the
the distance functions have an extremum at the same time.
So we attempt to find those $x \in \R^n$ for which there
are $s_{i}$ on the $M_i$ such that the  following
equations hold true:
\begin{gather}
 F_i(x,s) = \rVert x  -\gamma_i(s_i)\lVert \qquad F_i=F_j \qquad
1 \leq i,j \leq l 
 \nonumber \\
\label{eq:1}
\frac{\partial F_i}{\partial s_i} = 0 \in \R^{n-1} \qquad 1 \leq i \leq l
\end{gather}
Here $\gamma_i\colon M_i \mapsto \R^n$ are again embeddings for the
hypersurfaces $M_i$.
\begin{defi}
The $x_0$ for which such $\{ s_{i} \}_{1 \leq i \leq l}$ exist
make up the conflict set. 
\end{defi}
We measure not so much symmetry but more that what is in the middle.
\begin{figure}%[htbp]
\begin{center}
\includegraphics[height=5cm]{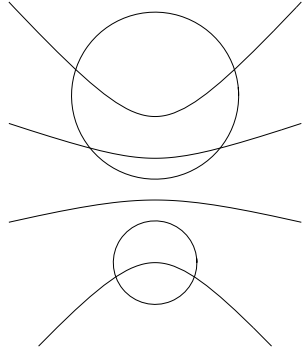}
\end{center}
\caption{A conflict set of two circles}
\end{figure}
The conflict set  of two lines is another pair of lines: nl.\ their
bisectors. For a line and a circle the conflict
set is already a more complicated object, 
while the conflict set of two circles generally will consist of
four conic sections. This is due to the fact that the distance function
from a fixed point to a circle mostly has two extrema. With other curves
such phenomena happen too. 
\newline
To be able to single out a certain component we define an 
oriented conflict set. Let the $M_i$ be oriented manifolds. Then
the orientation defines a direction for a flow induced by the
unit normal vectorfield $n_i \in NM_i$. 
\begin{defi}
The union of the intersection of the wavefronts at all times 
$t \in \R$ is called the oriented conflict set.
\end{defi}
The main result of this paper comes in two parts.
\begin{figure}[htbp]
\centering
\includegraphics[height=4cm]{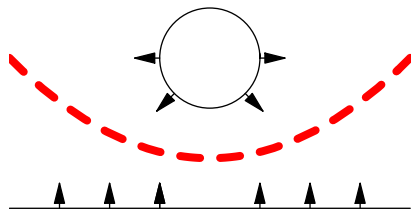}
\caption{An oriented conflict set of a line and a circle. The arrows indicate the orientations.}
\end{figure}
\begin{itemize}
\item The conflict set of $\{ M_i \}_{1\leq i \leq l}$
is for generic embeddings of the
surfaces $M_i$ a Legendrian manifold, that it is the 
projection of a smooth set in $PT^*\R^n$. 
\item The conflict set generically has the singularities that occur
in $n-l$ parameter families of 
wavefront in $\R^{n-l+2}$. If $n-l+2 \leq 6$ then 
the conflict set generically only has singularities that are combinations
of the well-known ADE singularities.
\end{itemize}
Our approach is such that it will entail several other results. Among
these are results concerning the center symmetry set, and others 
concerning the Gauss map. Also when $n=l$ we will define a sort of dual to
the conflict set, akin to the dual of a curve in projective space.
\newline
The paper has the following set up. In the first section we redefine
the conflict set as the intersection of so called big fronts. 
In the second section we use canonical relations to carry out
a number geometrical constructions. In the third section we
use the Thom transversality theorem to proof the first part
of out main result.  In the fourth section we recall some
results on Lagrangian and Legendrian singularities.
These mainly concern so-called phase functions. 
We use these to proof the second
part of our main result. 
\newline
The main results of this paper have been the subject of talks
held by the author at several conferences in the years 2000-2002.
\section{Statement of the main result}
\begin{figure}[hbtp]
\centering
\includegraphics[height=4cm]{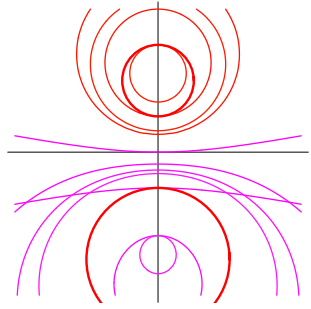}
\caption{Conflict sets for different values of $\eta$ }
\label{fig:6}
\end{figure}

Denote by $T^*\R^n\setminus 0 $ the slit cotangent bundle,
that is the cotangent bundle without the zero section.
Coordinates for the slit cotangent bundle  are 
$( x, \xi )$. The $\xi$ are coordinates in the fiber. 
Let $H_i\colon T^*\R^n\setminus 0 \rightarrow \R_{\geq 0}$
be $\Co{\infty}$ functions positively homogeneous of degree $1$,
and independent of $x$. If also the matrices
\begin{equation*}
  \frac{\partial^2 H_i^2(\xi)}{\partial \xi^2} 
\end{equation*}
are positive definite
these functions define a Finsler metric and
as Hamiltonians they define trajectories. 
These trajectories are particularly simple. They are
straight lines. More precisely we have that the
time that it takes to travel from $p_0$ to $p_1$ in
$\R^n$ is a function whose squared value is a $\Co{\infty}$
function on $\R^n\times \R^n$ minus the diagonal and whose first derivative
is never zero.
\newline
To a certain extent we could drop the condition that 
the $H_i$ are translation invariant or the condition 
that the $H^2_i$ only take on positive values but that
would lead us too far afield.
\newline
Now let $M_i$ be smoothly embedded manifolds of
codimension $1$ in $\R^n$. We could take smoothly 
embedded submanifolds of any codimension, but we consider
wavefronts emanating from these submanifolds and after
some time they become submanifolds of codimension $1$.
\newline
Denote coordinates on 
\begin{equation*}
 T^*\R^{n+1}\setminus  0 
\end{equation*}
by $(x,\xi,t,\tau) $.
\newline
For each of the $M_i$ we can define the ( oriented ) big
wavefront, by means of a map 
\begin{eqnarray*}
\Psi \colon T^*(\R^{n+1} ) \setminus 0
& \rightarrow &  \colon T^*(\R^{n+1} ) \setminus 0 \\ 
\Psi \colon (x,\xi,t,\tau) & \rightarrow & (\exp(t X_H),t,\tau-H)
\end{eqnarray*}
\begin{defi}\label{sec:stat-main-result-1}
The image of
\begin{equation*}
\R\times\{0\} \times ( N^*M_i  \cap \{ H_i =1 \} )
\subset T^*\R \times T^*\R^n\setminus  0 
\end{equation*}
by the map $\Psi$  or rather
its image under multiplication in the fibers by $\lambda\in\R_{>0}$,
is called the big wavefront. 
It is  denoted by $N^*M_i^h$.
\end{defi}
If the $M_i$ are orientable  we can consider a section 
$n_i$ of $N^*M_i$ contained in $\{H_i=1\}$. 
\begin{defi}
The image of
\begin{equation*}
\R\times\{0\}  \times \Image(n_i)
\subset T^*\R \times T^*\R^n\setminus  0 
\end{equation*}
by the map $\Psi$  or rather
its image under multiplication in the fibers by $\lambda\in\R_{>0}$,
is called the oriented big wavefront. 
It is  denoted by $N^*M_i^b$.
\end{defi}
The ( oriented ) big wavefront are conic Lagrangian submanifolds
of 
\begin{equation*}
T^*( \R \times \R^n) \setminus \{ 0 \}
\end{equation*}
We can now define a more general conflict set in a Hamiltonian
context.
\begin{defi}
The oriented conflict set of $(M_i,n_i,H_i)$, $1 \leq i \leq l$
is 
\begin{equation*}
C\left( \left(M_i,n_i,H_i \right)_{i=1, \cdots , l} \right)=
\pi_n \left( \bigcap_{i=1}^l \pi_{n+1} \left( N^*M_i^b \right) \right)  
\end{equation*}
where $\pi_n(x,t)=x$ and $\pi_{n+1}(t,\tau,x,\xi) = (t,x)$.
\end{defi}
\begin{rem}
The conflict set of $\{ (M_i, H_i)\}_{1\leq i \leq l} $ is
defined similarly. Also the sets $\pi_{n+1}(N^*M_i^h)$ and
$\pi_{n+1}(N^*M_i^b)$ are sometimes called ``the graph of the 
time function'', see \cite{MR93b:58019}. The ``time function''
is of course not really a function. It is multi-valued.
\end{rem}
Examples of these concepts are readily provided. For instance
we can consider the Hamiltonians
\begin{equation*}
  H_1 = \lvert \xi \rvert ~~ H_2= \eta \lvert \xi \rvert ~~ \eta \in \R_{>0}
\end{equation*}
for two circles. If we draw the conflict set for different values
of $\eta$, we get a plethora of curves , see figure~\ref{fig:6}. 
\newline
Denote by $\Emb(M_i,\R^n)$ the space
of embeddings of the manifold $M_i$ in $\R^n$. This space is
an open subset of $\Co{\infty}(M_i, \R^n)$.
To take into
account the $l$ embeddings we deal with we introduce the
space 
\begin{equation*}
  \bigoplus_{i=1}^l \Emb(M_i,\R^n)
\end{equation*}
The main result of our paper reads as follows. 
\begin{ste}\label{sec:stat-main-result}
For a residual subset of  $\bigoplus_{i=1}^l \Emb(M_i,\R^n)$
the ( oriented ) conflict set can be realized as the projection
of a conic Lagrange submanifold in $T^*\R^n$ . In addition if
the $M_i$ are compact the generic singularities 
of the ( oriented ) conflict set are the generic singularities
of $n-l$ parameter families 
of fronts in $\R^{n-l+2}$ embedded in $\R^n$. In particular
if $2 \leq n-l+2 \leq 6 $ the ( oriented ) conflict set
generically only has combinations of simple singularities of $ADE$ type.
\end{ste}
\begin{rem} The dimension of the conflict set in $\R^n$ is $n-l+1$.
\end{rem}
We point out one particularly nice consequence of this
theorem. If we have $n$ surfaces in $\R^n$ the conflict set
is one dimensional and 
generically only has  self-intersections and cusps as singularities.
\section{Canonical relations and associated geometrical constructs}
We will need the notion of a canonical relation between two 
symplectic manifolds. 
\begin{defi}
A  canonical relation  between 
two symplectic manifolds $\{ ( S_i, \omega_i ) \}_{i=1,2}$
is a  Lagrangian submanifold of
\begin{equation*}
  S_1 \times S_2 , \pi_1^*\omega_1 - \pi_2^*\omega_2
\end{equation*}
The $\pi_i$ denote projections $S_1\times S_2 \rightarrow S_i$.
\end{defi}
Graphs of symplectomorphisms $\chi\colon S_1 \rightarrow S_2$ are
important examples of canonical relations.
Symplectomorphisms can be composed, canonical relation
can be composed as well. The composition of symplectomorphisms is
a special case of a theorem of H\"{o}rmander, see \cite{Hormander:ALDO},
chapter 21.
\begin{ste}\label{sec:canon-relat-assoc}
Let $S_i, i=1\cdots 3 $ be three symplectic manifolds. Let $G_1$ be
a canonical relation between $S_1$ and $S_2$ and $G_2$ one between
$S_2$ and $S_3$. If $G_1 \times G_2$ intersects 
$S_1 \times \Delta (S_2) \times S_3$ transversally then the image $G_3$
under the projection $S_1 \times S_2  \times S_2  \times S_3  
\mapsto S_1  \times S_3 $ is a canonical relation between $S_1$ and $S_3$.
We call it the composition $G_1\circ G_2$ of $G_1$ and $G_2$.
\end{ste}
This theorem allows one to create new Lagrangian manifolds
from old. Canonical relations also formalize remarks of
Arnol'd who has regularly stated that in symplectic geometry one 
``decreases dimensions'' by ``sectioning and projection'', see for
instance \cite{MR93b:58019}. 
\newline
Often we will be interested in a particular case where
$S_3$ is a point. So we use the next proposition, that rephrases 
the demand in the theorem of H\"{o}rmander. Note that a 
canonical relation between $S_3$ and a point is just a  Lagrange 
manifold in $S_3$.
\begin{prop}\label{sec:canon-relat-assoc-1}
Let $G_1$ be a  canonical relation between $S_1$ and $S_2$ whose projection to
$S_2$ is an immersion and let $G_2$ be a  canonical relation between $S_2$ and a point. Then the
composition $G_1 \circ G_2$ is a  canonical relation , and thus a Lagrangian manifold in $S_1$, if 
$\pi_2(G_1) \pitchfork G_2$. 
\end{prop}
\begin{proof}
We need that 
\begin{equation}\label{eq:2}
G_1 \times G_2 \pitchfork S_1 \times \Delta(S_2)
\end{equation}
This intersection is contained in the graph of the 
projection $\pi_2\colon G_1 \rightarrow S_2$. We have 
that (\ref{eq:2}) holds iff.\
\begin{equation*}
  \mathrm{gr}(\pi_2) \pitchfork G_1 \times G_2
\end{equation*}
this in turn is true iff.\ 
\begin{equation*}
  \pi_2(G_1) \pitchfork G_2
\end{equation*}
\end{proof}
Now the theorem of H\"{o}rmander has a more familiar interpretation:
through this manifold $\pi_2( G_1)$ we can pull back 
Lagrangian manifolds from $S_2$ to $S_1$.
\subsection{Conflict sets}
Our foremost example of this procedure will 
be where
\begin{equation*}
  S_1 = T^*\R^{n+1} \setminus  0 
\end{equation*}
As coordinates for $S_1$ we will take 
\begin{equation*}
\left( \bar{x} , \bar{\xi} \right) 
~~\bar{x} = (x_0,x)  \in \R^{1+n} ~~ \bar{\xi} = ( \xi_0, \xi) \in \R^{1+n} 
\end{equation*}
or $(\bar{y}, \bar{\eta})$. In $T^*\R^n$ we will write 
$(x,\xi)$ or $(y,\eta)$.
\newline
If we put $S_2 = ( S_1 )^l$ then conflict sets can be
constructed by means of a canonical relation as follows. 
Set $G_1 \subset S_1 \times S_2$ 
\begin{equation*}
\left\{ \bar{y}, \bar{\eta} ,  \bar{x}_1, \bar{\xi}_1 ,\cdots ,
\bar{x}_l, \bar{\xi}_l \mid  \bar{\eta}= \sum_{i=1}^l
 \bar{\xi}_i , ~~\bar{y}=\bar{x}_i , ~~ 1 \leq i \leq l \right\}
\end{equation*}
The manifold $G_1$ is clearly conic Lagrange and also it projects
as an embedding to $S_2$. For $G_2$ we take the product of the
big wavefronts 
\begin{equation*}
  G_2=\times_{i=1}^l  N^*M_i^h \textrm{ or }  G_2=\times_{i=1}^l  N^*M_i^b
\end{equation*}
The manifold $\pi_2(G_1)$ is the diagonal in $\left( \R^{n+1}\right)^l$
together with all cotangent vectors
\begin{equation*}
\pi_2(G_1) = T^*_\Delta ( S_1^l )
\end{equation*}
The lifted conflict set $L^h$ ( or $L^b$ in the
oriented case )  is the pull-back to $S_1$ by 
$ T^*_\Delta ( \R^{(n+1)l} )$. The projection of $L_h$ to
the base $\R^{1+n}$  is the graph of the time function
on the conflict set.
\begin{ste}\label{sec:conflict-sets}
$L^h$ is conic Lagrange if 
\begin{equation}
\label{eq:3}
\times_{i=1}^l  N^*M_i^h \pitchfork T^*_\Delta ( \R^{(n+1)l} )
\end{equation}
\end{ste}
\begin{proof} Clear from the above remarks. \end{proof}
\begin{rem}\label{maxrankremark}
The criterion in equation (\ref{eq:3}) is quite
computable in practice. Suppose that the time functions for
each of the big wave fronts take on the form 
\begin{equation*}
t=F_i(x,s_i) , ~~ \frac{\partial F_i}{\partial s_i} =0, ~~ 1 \leq i \leq l
\end{equation*}
then the intersection in equation (\ref{eq:3}) is transversal iff.\
\begin{equation*}
F(x,\lambda,s_1 , \cdots , s_l) =\sum_{i=1}^{l-1}
\lambda_{i}(F_i(x,s_i)- F_{i+1}(x,s_{i+1}) )
\end{equation*}
is a phase function, i.e. the matrix 
\begin{equation*}
  \diffd_{s,\lambda,x}( F, \diffd_{s,\lambda} F )
\end{equation*}
has maximal rank there where 
\begin{equation*}
( F, \diffd_{s,\lambda} F ) = 0
\end{equation*}
In case the Hamiltonians are just Euclidean metrics this leads
to a lot of computable examples, in the spirit of \cite{Porteous:1994}. 
For instance, it is verified that with 3 surfaces in $\R^3$ only one of the
wavefronts can have a cuspidal edge, if this maximal rank criterion is to
hold. All more degenerate cases lead to,
albeit interesting, examples of non-Legendrian behavior.
\end{rem}
We now have a conflict set $L^h$ in $S_1$ and this object is
the most important one. First,
the projection of $L^h$ to $\R^{n+1}$ has the same 
Legendrian singularities as the projection of 
the corresponding object in $T^*\R^n\setminus\{ 0\} $ to
$\R^n$, and second, there is associated to the conflict
set a ``kite curve'' which can be constructed directly from
$L^h$. The kite curve will be treated below in paragraph \ref{sec:kite-curve}.
\newline
To pull $L^h$ into $T^*\R^n\setminus\{ 0\} $ we need to apply the
``sectioning and projection''. The section is $\eta_0=0$ and 
the projection is along the $y_0$ axis. The $y_0$ is axis is the
time axis and as we saw that $\pi_{n+1}(L^h)$ is the graph of the
time function on the conflict set it is not surprising that this
projection is immersive and that it induces no extra singularities.
The conic canonical relation we use is
\begin{equation}\label{eq:26}
  G_1 = 
\{ (x,\xi,\bar{y},\bar{\eta}) \mid 
x=y,~~ \xi = \eta, \eta_0=0 \} \subset T^*\R^n\setminus\{ 0\} \times S_1 \textrm{ and } G_2=L^h
\end{equation}
Using proposition \ref{sec:canon-relat-assoc-1} we obtain that if
\begin{equation}
  L^h \pitchfork W \textrm{ where }
W= \{ (\bar{y},\bar{\eta} ) \mid \eta_0=0 \} 
\end{equation}
we can pull back $L^h$ to  $T^*\R^n\setminus\{ 0\} $.
\begin{lem}
$\times_{i=1}^lN^*M_i^h
\pitchfork T_\Delta^*\left( S_1^l \right)
\Rightarrow W \pitchfork L^h$
\end{lem}
\begin{proof} Suppose that we did not have $W\pitchfork L^h$.
Because $W$ is a hypersurface that would mean that 
at some point $p$ in $L^h$ the tangent space $TL^h$ would 
be contained in $TW$. So it would hold 
\begin{equation*}
\langle (0,0,0,\delta \eta_0) , \vec{w} \rangle =0 , ~~ \forall \vec{w}\in TL^h
\end{equation*}
and consequently
\begin{equation*}
  \omega( (0,0,\delta y_0,0),\vec{w} )=0  , ~~ \forall \vec{w}\in T_pL^h
\end{equation*}
with $\omega$ being the canonical symplectic structure.
But $L^h$ is Lagrangian, so we'd have that
this vector $(0,0,\delta y_0,0 ) \in TL^h$.
But that is clearly impossible. 
\end{proof}
\begin{cor}
The conflict set is conic Lagrange if (\ref{eq:3}) holds.
\end{cor}
Our main theorem says that not only generically the conflict set
is a Legendrian manifold, but also that its singularities are singularities
of fronts in $\R^{n-l+2}$. So we want to pull back not
just to $\R^n$ but to an $n-l+2$ dimensional
subspace of $\R^n$. This can be done in more or less the same 
way as in the above where we projected along along the ``time'' axis.
front in an $n-l+2$ dimensional space. 
\newline
We can project in a Legendrian way along some direction $v$ if 
\begin{equation*}
  W(v) = \{ (x , \xi ) \in T^*\R^n 
\mid \langle v , \xi \rangle = 0 \qquad \lVert \xi \rVert = 1 \}
\end{equation*}
intersects transversely with $L^h$. This was done in the above along
the time axis.
\newline
More generally we can section and project along a subspace $V$ spanned
by $\{ v_1 , \cdots , v_k \}$, if 
\begin{equation*}
  W(V) = \{  (x , \xi ) \in T^*\R^n \mid
\langle v_i , \xi \rangle = 0 \quad
i=1 , \cdots , l \quad \lVert \xi \rVert = 1 \}
\end{equation*}
intersects $L^h$ transversely. The assertion is proven
by using a canonical relation as the previous $G_1$. In
$\R^n$ the conflict set is $n-l+1$ dimensional. If it is 
Legendrian ``the fiber'' has dimension  $l-1$. A maximum of
$l-2$ directions can thus additionally be sectioned 
away. We end up in an $n-l+2$ dimensional space, 
as stated in the main theorem.
\subsection{The kite curve}\label{sec:kite-curve}
Associated to the conflict set is the kite curve. It is some sort of
dual to the conflict set. This occurs
when $l=n$.
Suppose we are at a point of the conflict set
that corresponds to $l$ Morse extrema of the
time function and such that the normals to
the front are affinely in general position.
\begin{figure}[htbp]
\centering
\includegraphics[height=4cm]{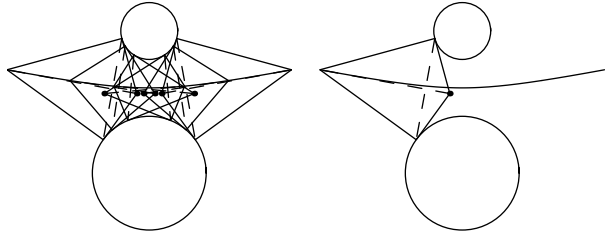}
\caption{Oriented conflict sets and some kites}\label{fig:8}
\end{figure}
The conflict set is smooth there and at the
$l$ basepoints on the $M_i$ we have tangent
planes, that are naturally thought of as affine
hyperplanes in $\R^n$. They usually intersect
in a point. This point traces out a curve in
$\R^n$.
\newline
The kite curve can most conveniently be thought
of the intersection of the tangent developable
of $\pi_{n+1}(L^h)$ with a plane $t=\mathrm{constant}$,
see figure \ref{fig:kiteconstr}. 
This explains immediately that the kite curve
is a line whenever the $M_i$ are all spheres
( in sufficiently general position ).
\newline
To form a good idea of the kite curve one could
consult \cite{Siersma1}. There the case with 2
curves in $\R^2$ is thoroughly examined. The
kite curve is exhibited there as the locus of
points where two tangent lines have equal length.
As said for two circles in the plane the kite
curve is a straight line, see the pictures in
\ref{fig:8}.
\newline
If one wants to, the kite curve can also be
constructed at singular points of the conflict set.
First of all we are interested in the intersection
of the big wave fronts, so we need:
\begin{equation}
  \label{eq:4}
   \bar{x}_i=\bar{x}_j , ~~ 1 \leq i , j \leq l
\end{equation}
Our $y$ variable should be such that it is in
the intersection of the tangent planes to the big
fronts. So the vector that runs from the point
$(y,0)$ to $\bar{x}_i$ should lie in all the 
tangent planes. So the line from $(y,0)$ to
$\bar{x}_i$ should be orthogonal to each
$\bar{\xi}_i$. That is:
\begin{equation}
  \label{eq:6}
  \langle x -y , \xi_i \rangle +  x_0  \xi_{i,0}  =0 
  , ~~ 1 \leq i \leq l =n
\end{equation}
Equations (\ref{eq:4}) and (\ref{eq:6}) define a 
set in $\R^n$. Possibly this is a
Legendrian curve, but this is not so clear. 
Further on we will consider a more natural candidate 
for a ``kite curve''. 
\newsavebox{\mygraphic}
\sbox{\mygraphic}{\includegraphics[width=9cm]{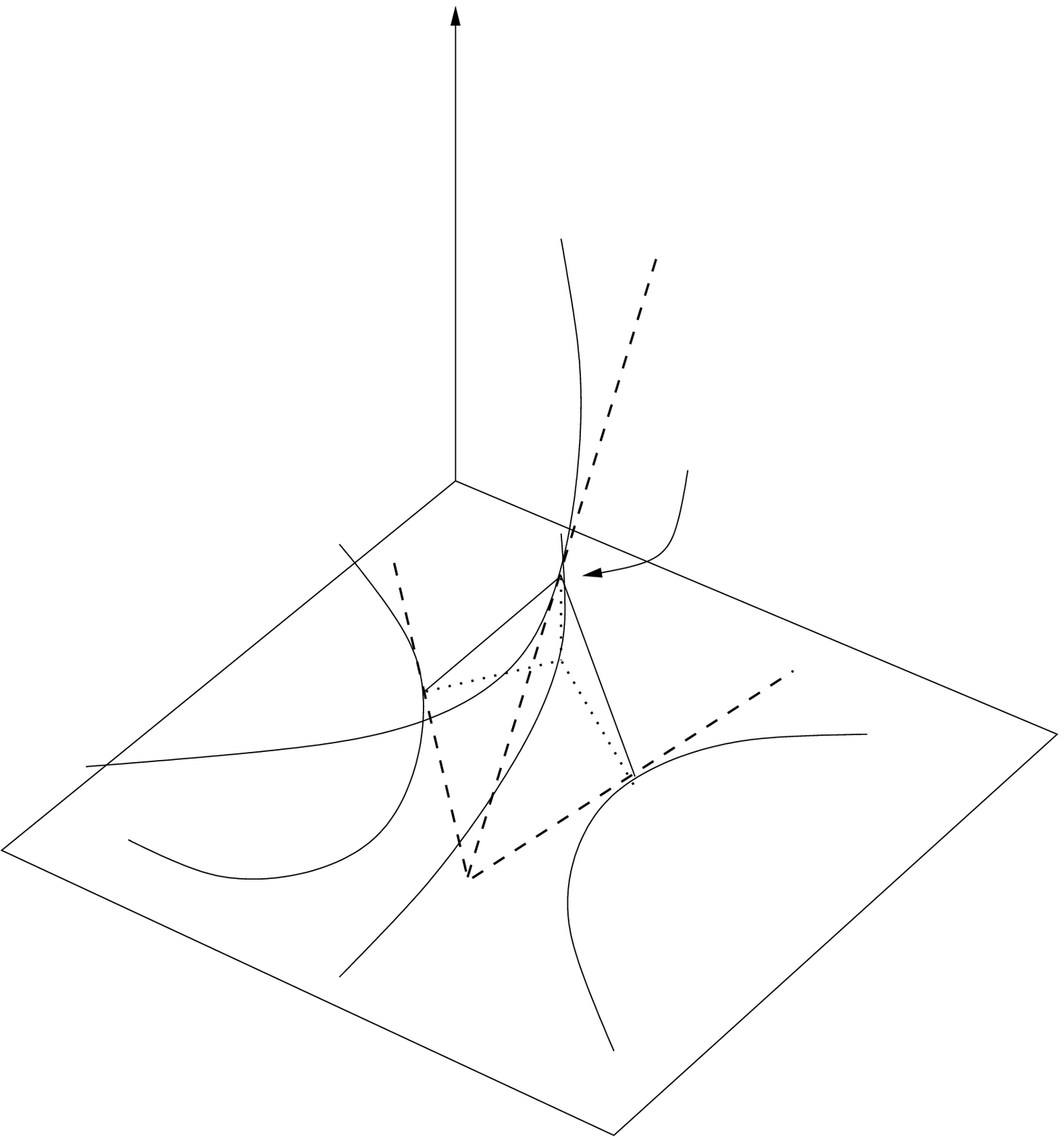}}
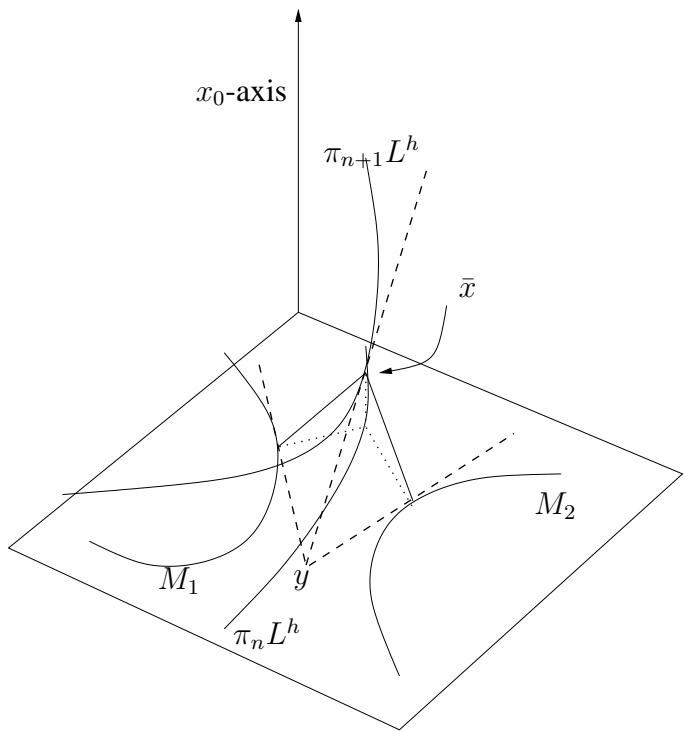
\begin{figure}[hbtp]
\centering
\setlength{\unitlength}{1.0cm}
\begin{picture}(9,10)
\put(0,0){\usebox{\mygraphic}}
\put(2.5,8.5){\makebox(0,0)[l]{$x_0$-axis}}
\put(6,5.9){\makebox(0,0)[l]{$\bar{x}$}}
\put(3.8,2.0){\makebox(0,0)[l]{$y$}}
\put(2,2){\makebox(0,0)[l]{$M_1$}}
\put(7,3){\makebox(0,0)[l]{$M_2$}}
\put(4.2,7.7){\makebox(0,0)[l]{$\pi_{n+1}L^h$}}
\put(3,1.3){\makebox(0,0)[l]{$\pi_nL^h$}}
\end{picture}
\caption{The construction of the kite curve.}\label{fig:kiteconstr}
\end{figure}
\subsection{Gauss maps and parallel tangent planes}
The Gauss map can also be used 
to define sets that measure symmetry. 
The center symmetry set is defined as the locus of the
midpoints of chords connecting two distinct points with
parallel tangent space on a surface 
$M_1 \subset \R^n$. Apparently this was introduced by Giblin,
see \cite{MR2001e:58054}.
The center symmetry set of a conic section is a point. 
\newline
Clearly we can repeat this with two distinct manifolds, looking
for pairs of parallel tangent planes. As set we can consider
\begin{enumerate}
\item the midpoints of the chords connecting the tangent planes
\item the normals themselves as a subset of the family of oriented
lines in $\R^n$, $T^*S^{n-1}$
\end{enumerate}
\begin{defi} The midpoints of the chords form the \textbf{center set}
and the chords themselves form the \textbf{normal chord set}.
\end{defi}
These constructions can both be carried out with 
canonical relations. The first curve will turn out to be 
Legendrian iff.\ the second is Lagrangian.
\newline
We do a recap of the construction of the space of oriented lines in
$\R^n$. A directed line in $\R^n$ has a direction, that is a unit
vector $v$ in $S^{n-1}$. At $v\in S^{n-1}$ there is a tangent plane. This 
tangent plane can be identified with a plane in $T_0\R^n$.
The intersection point of this tangent plane with the directed line
leaves us with a vector in $T_vS^{n-1}$. On the tangent space we
have the Legendre transform that maps the vector in $T_vS^{n-1}$
to $T^*_vS^{n-1}$.
\newline
On a hypersurface $M$ in $\R^n$ we have the 
Gauss map. We can view it as a map to
$T^*S^{n-1}$: assign to $p\in M$ the normal as a directed line
in $\R^n$. It is a standard theorem that the image of the
Gauss map is Lagrangian, see again \cite{MR93b:58019}.
\newline
One can proof this theorem by proving that the
set
\begin{gather}
v, \mu , x, \xi \in T^*S^{n-1}\times T^*\R^n \textrm{ such that } \nonumber\\
\label{eq:17}
v = \frac{\xi}{\lVert \xi \rVert}, ~~ \mu = x -
\frac{ \langle x ,\xi \rangle \xi}{\lVert \xi \rVert^2} , ~~
\lVert \xi \rVert = 1
\end{gather}
is a canonical relation between $T^*S^{n-1}$ and $T^*\R^n$.
The pull-back to $T^*S^{n-1}$ of the conormal bundle to
a manifold is the image of the Gauss map. We could try to 
pull back any Lagrangian manifold to $T^*S^{n-1}$ by computing
in this way ``the image of the Gauss map''. Using
proposition \ref{sec:canon-relat-assoc-1} we see that 
this only makes sense if the Lagrange manifold is $T^*\R^n$ 
lies transverse to all the level sets of $\lVert \xi \rVert$ ,
that is, ``the image of the Gauss map'' is defined for
all conic Lagrange manifolds. This is a rather important
remark, we formulate it in a theorem - for which we claim
no originality whatsoever:
\begin{ste}
For each Legendre manifold in $\R^n$ we can define
an ``image of the Gauss map''. If the Legendre manifold
is the conormal bundle of a smooth submanifold of
codimension $1$ in $\R^n$ this coincides with the
usual image of the usual Gauss map.
\end{ste}
\subsection{The center set}
Let us now look at the center set. We introduce a 
canonical relation between
\begin{equation*}
  S_1=T^*\R^n\setminus 0 
\end{equation*}
and 
\begin{equation*}
  S_2= (S_1)^2
\end{equation*}
We will as in the above use coordinates 
\begin{equation*}
  (y , \eta, x_1, \xi_1 , x_2, \xi_2 )
\end{equation*}
Define $G_1$ by 
\begin{equation}
  \label{eq:11}
  y = \frac{x_1+x_2}{2} , ~~ \eta=2\xi_1, \xi_1=\xi_2
\end{equation}
The manifold $G_1$ is clearly a conic canonical relation and $G_1$ projects immersively
into $S_2$. For $G_2$ we take the product of the conormal
bundles of $M_1$ and $M_2$. The proposition \ref{sec:canon-relat-assoc-1}
can be applied to yield that 
\begin{ste}\label{sec:centre-set}
The center set is conic Lagrange if 
\begin{equation}
\label{eq:7}
N^*M_1 \times N^*M_2 \pitchfork \{(x_1,\xi , x_2 , \xi \}
\end{equation}
\end{ste}
\begin{rem}
Note the reciprocity between (\ref{eq:7}) and (\ref{eq:3}). 
The criterion for the conflict set deals with a diagonal in the
base and the criterion for the center set deals with a diagonal in the
fiber
\end{rem}
\begin{rem}
Clearly, lots of other interesting and less interesting
sets can be constructed in this way. We could take $l$ manifolds $M_l$ and
consider the relation
\begin{equation}\label{eq:25}
y=\sum_{i=1}^l a_i x_i  \text{\phantom{aap}} 
\eta a_i = \xi_i
\end{equation}
where the $a_i$ are a set of nonzero numbers. This again is
a canonical relation. If the product of the conormal bundles of
the $M_i$ is transverse to the ``diagonal in the fiber'' as in
(\ref{eq:7}) then the resulting set is a Legendre manifold.
For instance with three surfaces we could take the centroid of a triangle.
\end{rem}
\subsection{The normal chord set}
Next comes the normal chord set. Denote $\nu_i$ the Gauss map 
from $N^*M_i$ to $T^*S^{n-1}$. We have an image of
the product of the two Gauss maps, the chords we 
use are on the diagonal. After the many examples 
above the following theorem is obvious:
\begin{ste}
If 
\begin{equation}
  \label{eq:12}
 \nu_1(N^*M_1) \times \nu_2(N^*M_2) \pitchfork \{ (v,\mu_1, v , \mu_2 )\} 
\end{equation}
the normal chord set is Lagrangian.
\end{ste}
\begin{rem}
For the normal chord set we can write down maximal rank criteria as
we did for the conflict set in remark \ref{maxrankremark}.
A phase function for the image of the Gauss map is
\begin{gather}
  F_i\colon S^{n-1}\times M_i \rightarrow \R \nonumber \\
 (v, s ) \mapsto \langle v , \gamma(s) \rangle
  \label{eq:13}
\end{gather}
The image of the Gauss map is described by 
\begin{equation}\label{eq:14}
  \frac{\partial F_i}{\partial s} =0 
\end{equation}
To get a maximal rank criterion under which the normal chord set
is Lagrangian we use as in the remark \ref{maxrankremark}
a special phase function:
\begin{equation*}
   F_1(v,s_1) +  F_2 (v, s_2 )
\end{equation*}
And the maximal rank criterion that is equivalent to
the transversality in (\ref{eq:12}) is that the matrix
\begin{equation*}
\diffd_{ v, s_1 ,s_2 }
\left( \diffd_{ s_1, s_2} F \right) 
\end{equation*}
has maximal rank there where 
\begin{equation*}
   \diffd_{s_1, s_2} F = 0
\end{equation*}
If so, the normal chord set is Lagrangian. If one chooses local coordinates
on $S^{n-1}$, as is done in \cite{MR84a:58001}, this is a nicely
computable criterion. 
\end{rem}
\subsection{The kite curve revisited}
As a definition of the kite curve it seems more useful to consider
the image in $T^*S^n$ of the 
lifted conflict set $L^h$: the image of its Gauss map. This is also defined
when $l < n$. The kite curve can be reconstructed from it, when $l=n$.
\newline
If we define the kite curve in this way we can summarize our reasoning 
in a tentative diagram, that illustrates the dualities mentioned.
\newline
\setlength{\unitlength}{3000sp}
\begingroup\makeatletter\ifx\SetFigFont\undefined%
\gdef\SetFigFont#1#2#3#4#5{%
  \reset@font\fontsize{#1}{#2pt}%
  \fontfamily{#3}\fontseries{#4}\fontshape{#5}%
  \selectfont}%
\fi\endgroup%
\begin{picture}(8400,4000)(0,-4500)
\thinlines
\put(2401,-961){\line( 0,-1){3600}}
\put(4801,-961){\line( 0,-1){3600}}
\put(7201,-961){\line( 0,-1){3600}}
\put(601,-2161){\line( 1, 0){8400}}
\put(601,-3361){\line( 1, 0){8400}}
\thicklines
\put(5801,-3061){\vector( 0,-1){600}}
\put(3801,-3061){\vector( 0,-1){600}}
\put(3001,-1561){\makebox(0,0)[lb]{\smash{\SetFigFont{12}{14.4}{\rmdefault}{\mddefault}{\updefault}Diagonal in the}}}
\put(3001,-1801){\makebox(0,0)[lb]{\smash{\SetFigFont{12}{14.4}{\rmdefault}{\mddefault}{\updefault}base}}}
\put(5401,-1561){\makebox(0,0)[lb]{\smash{\SetFigFont{12}{14.4}{\rmdefault}{\mddefault}{\updefault}Diagonal in the}}}
\put(5401,-1801){\makebox(0,0)[lb]{\smash{\SetFigFont{12}{14.4}{\rmdefault}{\mddefault}{\updefault}fibre}}}
\put(5701,-2761){\makebox(0,0)[lb]{\smash{\SetFigFont{12}{14.4}{\rmdefault}{\mddefault}{\updefault}Centre set}}}
\put(5701,-3961){\makebox(0,0)[lb]{\smash{\SetFigFont{12}{14.4}{\rmdefault}{\mddefault}{\updefault}Chord set}}}
\put(3001,-2761){\makebox(0,0)[lb]{\smash{\SetFigFont{12}{14.4}{\rmdefault}{\mddefault}{\updefault}Conflict set}}}
\put(3001,-3961){\makebox(0,0)[lb]{\smash{\SetFigFont{12}{14.4}{\rmdefault}{\mddefault}{\updefault}Kite curve}}}
\put(901,-2761){\makebox(0,0)[lb]{\smash{\SetFigFont{12}{14.4}{\rmdefault}{\mddefault}{\updefault}Legendrian}}}
\put(901,-3001){\makebox(0,0)[lb]{\smash{\SetFigFont{12}{14.4}{\rmdefault}{\mddefault}{\updefault}lift is smooth}}}
\put(901,-3961){\makebox(0,0)[lb]{\smash{\SetFigFont{12}{14.4}{\rmdefault}{\mddefault}{\updefault}Lagrangian}}}
\put(901,-4201){\makebox(0,0)[lb]{\smash{\SetFigFont{12}{14.4}{\rmdefault}{\mddefault}{\updefault}lift is smooth}}}
\put(7801,-2761){\makebox(0,0)[lb]{\smash{\SetFigFont{12}{14.4}{\rmdefault}{\mddefault}{\updefault}in $T^*\R^{n+1}$}}}
\put(7801,-3001){\makebox(0,0)[lb]{\smash{\SetFigFont{12}{14.4}{\rmdefault}{\mddefault}{\updefault}(resp.\ $T^*\R^n$)}}}
\put(7801,-3961){\makebox(0,0)[lb]{\smash{\SetFigFont{12}{14.4}{\rmdefault}{\mddefault}{\updefault}in $T^*S^n$ }}}
\put(7801,-4201){\makebox(0,0)[lb]{\smash{\SetFigFont{12}{14.4}{\rmdefault}{\mddefault}{\updefault}(resp.\ $T^*S^{n-1}$)}}}
\end{picture}
\section{Genericity of transversality conditions}
We proof that for a residual set of embeddings in 
  $\bigoplus_{i=1}^l \Emb(M_i,\R^n)$ the transversality
conditions that ensure that the conflict set and the 
center set ( resp.\ the kite curve and the normal chord set )
are Legendrian ( resp.\ Lagrangian ), are satisfied.
\subsection{That the conflict set is generically Legendre}
To proof the genericity of the criterion (\ref{eq:3}) we 
will make use of the map that defines the big wavefront,
described in definition \ref{sec:stat-main-result-1}.
This is a map 
\begin{equation}\label{eq:18}
\times_{i=1}^l (M_i \times (\R\setminus 0) \times \R )
\rightarrow T^*(\R^{n+1})^l
\end{equation}
More precisely, we associate such a map to
each $(\gamma_1, \cdots , \gamma_l)$.
\begin{equation}\label{eq:19}
\bigoplus_{i=1}^l \Emb(M_i,\R^n) \rightarrow 
\Co{\infty}(
\times_{i=1}^l (M_i \times (\R\setminus 0) \times \R), 
 T^*(\R^{n+1})^l)
\end{equation}
To simplify matters we look at each embedding individually. 
\begin{equation}\label{eq:24}
\Emb(M_i,\R^n) \rightarrow 
\Co{\infty}(
M_i \times (\R \setminus 0 ) \times \R , T^*(\R^{n+1}) )
\end{equation}
These can be put in a family. Namely just translate the
embeddings by a ( small ) vector $e_i$. For simplicity
drop the index $i$.
\begin{gather}
\R^n \times \Emb(M,\R^n) \rightarrow 
\Co{\infty}(
M \times (\R \setminus 0 ) \times \R,
T^*(\R^{n+1})^l ) \nonumber \\
\left(e,  \gamma\colon M \rightarrow \R^n \right)
\longmapsto 
 e,s,\lambda , x_0 
\overset{\Psi}{\rightarrow }
\begin{pmatrix}
 \pi_x( \exp(x_0 X_H)( \gamma(s)  + e) ) \\ 
 \lambda \pi_\xi( \exp(x_0 X_H)( \gamma (s) + e) ) \\
x_0 \\  \lambda 
\end{pmatrix}
\end{gather}
We need that $l$ copies of maps $\Psi$ ,each for
a different $M_i$ map transversal to
$T^*_\Delta(\R^{1+n})^l$. We see that it suffices to
proof that $l$ copies of 
\begin{equation}\label{eq:20}
e,s,x_0, \lambda \rightarrow \pi_x( \exp(x_0 X_H)( \gamma(s)  + e) ) , x_0
\end{equation}
map transversal to the diagonal $\Delta \subset  (\R^{1+n})^l$.
Because 
\begin{equation*}
 \pi_x( \exp(x_0 X_H)( \gamma(s)  + e) ) =
  e + \pi_x( \exp(x_0 X_H)( \gamma(s)))
\end{equation*}
this is clear: the derivatives for the $e$ vectors and the time 
variable $x_0$ already cause the maximal rank to be attained.
So the product of $l$ maps (\ref{eq:20}) is submersive 
onto $(\R^{n+1})^l$.  We have thus shown that 
\begin{equation}\label{eq:21}
( \R^n )^l \times  \bigoplus_{i=1}^l \Emb(M_i,\R^n) \rightarrow 
\Co{\infty}(
\times_{i=1}^l (M_i \times (\R\setminus 0) \times \R), 
 T^*(\R^{n+1})^l)
\end{equation}
gives for all 
\begin{equation*}
 \left(  \gamma_1 , \cdots , \gamma_l  \right)
\end{equation*}
a family of mappings , parametrized by $(\R^n)^l$, whose
image is transverse to the closed manifold
$T_\Delta^*(\R^{1+n})^l$, or any closed submanifold of 
$T^*(\R^{1+n})^l$ for that matter. Thus most members of
this family are transversal to $T_\Delta^*(\R^{1+n})^l$.
It now follows from results of Abraham, see \cite{Wall},
that for a residual set of embeddings the transversality
condition is satisfied, which is exactly what we needed. 
For clarity we cite the theorem:
\begin{ste}[Thom Transversality theorem]\label{sec:that-conflict-set}
Let $A$ be a manifold of mappings , let $X$, $Y$ be 
manifolds. Let 
\begin{equation*}
\alpha \colon A \rightarrow \Co{\infty}( X, Y)
\end{equation*}
be a map such that 
\begin{equation}\label{eq:23}
\mathrm{ev}(\alpha) \colon A\times  X\times  Y
\end{equation}
is a smooth submersion at every $(a,x)\in A\times X$. 
Then for every closed submanifold $W\subset Y$ we have 
that 
\begin{equation*}
  \{ a \in A \mid ~~ \alpha(a) \pitchfork W \}
\end{equation*}
is a residual subset of $A$.
\end{ste}
In particular $\oplus_{i=1}^l\Emb(M_i, \R^n)$ is a manifold of
mappings, the map $\alpha$ we use is in
(\ref{eq:24}). For any fixed $(\gamma_1 , \cdots , \gamma_l)$ $l$ 
copies of the family in~(\ref{eq:20}) are submersive. 
These are all embeddings, so locally the map from 
(\ref{eq:23}) is a smooth submersion. 
\begin{prop}
The conflict set is generically Legendre and 
the kite curve in $T^*S^n$ is generically Lagrange.
\end{prop}
\begin{rem}
Roughly speaking the family of translations (\ref{eq:20}) produces all
first order perturbations. This family will be of much use
further on. 
\end{rem}
\subsection{That the center set is generically Legendre}
\begin{figure}
\centering
\includegraphics[height=4cm]{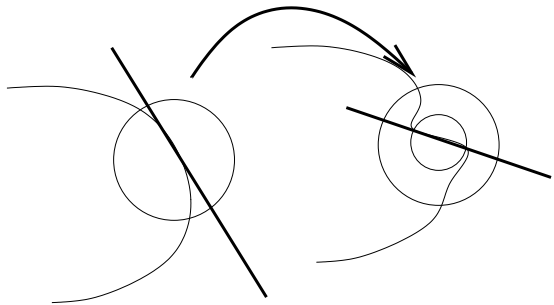}
\caption{The map $\phi_{r,A,p=\gamma(s^\prime)}\circ\gamma(s)$}\label{fig:turntangentspace}
\end{figure}

The aforementioned family seemingly can not be used
to proof that generically the center set is conic Lagrange.
We need a covering $\{U_\alpha\}$  of $M_1\times M_2$ and 
in each $\{ U_\alpha \}$ perturb the tangent space a
little, as indicated in figure \ref{fig:turntangentspace}.
It is enough to proof that , 
if $\vec{n}_i$ is the map that assigns the normal to
$M_i$, that the map $(\vec{n}_1 , \vec{n}_2)$ 
is transverse to the diagonal. We first show that
locally families exist that are indeed transverse to the diagonal.
Denote by 
\begin{equation*}
\phi_{r, A, p } , ~~ r \in \R , A \in SO(n,\R), p\in \R^n
\end{equation*}
a diffeomorphism, which is the identity on $\R^n$ where we're
outside the sphere of radius $2r$ round $p$ and equal to
$q\rightarrow A(x-q)$, inside a circle of radius $r$ 
round $q$. Now compose an embedding
$\gamma\colon M \rightarrow \R^n$ with the map
$ \phi_{r,A,p(\alpha)} $ and we get 
a map that in some environment $U_\alpha^\prime$  of
of $p(\alpha)\in M$ is submersive. Looking
at a product $ \phi_{r ,A
, q(\alpha)}\circ\gamma_1 ,\phi_{r^\prime, A^\prime ,p(\alpha)}
\circ \gamma_2$
we see that in a neighborhood $U_\alpha$ the 
transversality condition is satisfied.
Indeed, at  $p(\alpha)$ the normal looks
like $A\vec{n}$. 
\newline
One can pick a countable number of points $p(\alpha)$ 
such that the $U_\alpha$ cover $M_1\times M_2$ . We have proven:
\begin{ste}
For a countable intersection of open and dense
subsets of 
$$
\bigoplus_{i=1}^2 \Emb(M_i,\R^n)
$$
the center set is Legendrian ( and the normal chord
set therefore Lagrangian )
\end{ste}
\section{Proof of main theorem}
The arguments presented in this section are standard and so
we do not provide all the details. More details can be found in the 
excellent surveys \cite{Duistermaat:CPAM} and
\cite{Wall}. See also \cite{Arnold:SAD}.
\newline
The general idea of the proof of the theorem is that
we stratify for each $M_i$ some part $B_i$ 
of the graph of the time function $\pi_{n+1}(N^*M_i^h)$.
The strata  of $B_i$ correspond to singularity types 
of individual momental fronts. Then the intersections
of these graphs are generically such that they miss
the non-stratified part. So in the intersection
we will only meet singularities that are well-known
singularity types of individual fronts. 
\newline
Our proof will consist of purely local considerations, they can be
patched together as the transversality theorem.
\subsection{Stratification by codimension of the big wavefront}
We want to define equisingularity manifolds. They will be
defined using the notion of codimension for germs. We use
codimension wrt.\ to contact equivalence or  V-equivalence.
This is motivated by our Legendrian point of view. Two germs of
unfoldings are $V$-isomorphic iff.\ the germs of Legendrian
immersions they determine are equivalent, see \cite{Arnold:SAD}, 
\S 20.
\newline
We consider locally $x_0=F(x,s)$ the time function belonging
to the embedding $M\hookrightarrow \R^n$. That is we
consider it near $\bar{x}^\prime ,s^\prime \in \R^{n+1}$. Then we put 
$G(\bar{x},s)=x_0-F(x,s)$.
The equations
\begin{equation}
  \label{eq:15}
  \frac{\partial G}{\partial s} = 0, ~~ G=0
\end{equation}
define a ( germ of a ) surface $Z_1 \subset \R^{1+n+n-1}$,
near $\bar{x}^\prime ,s^\prime$.
We want to
consider closed parts of the surface $Z_1$, namely
those where the codimension of the germ 
\begin{equation}
  \label{eq:16}
  \jmath G \in \Co{\infty}(s^\prime), ~~
s \overset{\jmath G}{\rightarrow} G(\bar{x}^\prime,s)=x_0 
\end{equation}
has a suitable value:
\begin{equation*}
  1 \leq \codim \jmath G \leq N
\end{equation*}
We define the codimension by the dimension of the 
real vector space
\begin{equation*}
  \dim_\R
\frac{ \Co{\infty}( s^\prime ) }{ \left( G , \diffd_s G \right) }
\end{equation*}
As an example consider 
\begin{equation*}
  x_0 = s^3+x_1s+x_2 , ~~ \bar{x}^\prime=0 , ~~ s^\prime=0
\end{equation*}
The codimension to calculate is 
\begin{equation*}
  \dim_\R \frac{ \Co{\infty}(s^\prime) }{(s^3,s^2) } = 2
\end{equation*}
There are other points nearby that have codimension 2. These
are 
\begin{equation*}
\{x_0=a,~~ x_1=0,~~x_2=a, ~~s=0 \}
\end{equation*}
The projection of this manifold to the $\bar{x}$-space is
a affine space. This affine space has codimension 2 in $\R^{n+1}$. ( In the example
$n=2$ ) 
\newline
We would like this to hold in general. The manifold $Z_1$ should be stratified.
\newline
Put $N=\min(n-l+2,6)$. Then in the space of $N+2$-jets of germs
at $s^\prime$ in $M$, we have a part $\mathcal{C}$ that is stratified
according to the codimensions $1$ to $N$. The part of
$J^{N+2}(s^\prime)$ with 
codimension $>N$ is an algebraic variety that can be 
stratified in some canonical way.  
\newline
The surface $Z_1$ is to be divided into a 
part $\mathcal{B}$ and its complement $\complement\mathcal{B}$.
The complement should have codimension $>N$ and the 
part $\mathcal{B}$ should have a Whitney stratification 
such that on each stratum the codimension is 
constant.
\newline
To stratify $Z_1$ we consider the map
\begin{equation*}
   \bar{x},s  \overset{\jmath^{N+2}G}{\rightarrow}  \R^n \times J^{N+2}(s) 
\end{equation*}
If $\jmath^{N+2}G$ is transverse to the stratification of $\mathcal{C}$
and its complement then this induces the division of $Z_1$ in
$\mathcal{B}$ and $\complement\mathcal{B}$.
Each of the strata of $\mathcal{B}$ corresponds to finitely
many types of singularities and they project to $\R^{1+n}$ 
immersively. 
\newline
We next pass to intersections.
\newline
Above one $\bar{x}$ there might be
several pairs of $(\bar{x},s^{(i)}),~~i=1 \cdots r$.
They are only finitely many because $M$ is compact.
For a residual subset of $G$ and thus for a residual
subset of the embeddings $\gamma$ the projection 
$\pi^{r}\colon Z_1^{(r)} \rightarrow (\R^{1+n})^r$
is transverse to the diagonal stratification of $(\R^{n+1})^r$.
We can have $\pi^r(\mathcal{B})\pitchfork\mathcal{D}^{(r)}$
generically for $r=n$, which will be enough.
\newline
For this reason generically the stratification of $\mathcal{B}$ has
regular intersections relative $\pi\colon Z_1\rightarrow\R^{1+n}$.
\begin{lem}\label{sec:strat-codim-big}
For a residual set of embeddings $M\rightarrow\R^n$ the 
``graph of the time function'' $\pi_{n+1}(N^*M^h)$ has a
subset $\mathcal{B}$ of codimension $\leq N$ that is Whitney
stratified and whose strata correspond to singularity
types of individual momental fronts.
\end{lem}
\subsection{Intersecting $l$ big fronts}
We want the intersection of the $l$ bigfronts be such
that in the intersection we only find 
elements of $\bigcap_{i=1}^l\mathcal{B}_i$ 
and no points of one the complements $\complement\mathcal{B}_i$.
With the family~(\ref{eq:21}) we see that the
intersection of the strata of the bigfronts is transversal.
\newline
In fact the maximal codimension of a stratum of say
 $\pi_{n+1}(N^*M^h_1)$
that can appear in the intersection of the big fronts, appears
when the other $l-1$ big fronts are smooth, so this maximal
codimension is $\leq n+1 - (l-1)$. Thus if $N \leq n-l+2$
we have only strata of the $\mathcal{B}_i$ in the intersection.
This is where the condition $n-l+2\leq 6$ in the theorem 
comes in.
\newline
The last thing we need to know to hold generically
is alike what we needed to know for the 
projection $Z_1\rightarrow\R^{1+n}$, namely
that $\pi\colon\cap_{i=1}^l\pi_{n+1}(L^h)\rightarrow\R^n$
has regular intersections relative
$\pi_n\colon L^h\rightarrow \R^n$.
This is again achieved with the family~(\ref{eq:21}). 
\newline
If $n-l+2=7$ then it will be possible to find a stratum
of the complement of say $\mathcal{B}_2$ in the intersection
of the big fronts. Such a stratum can represent
a modulus. If we move $\pi_{n+1}(N^*M^h_2)$ a little the
stratum will still be there, but the singularity type
will have changed. We conclude that  $n-l+2 \leq 6$ are
the nice dimensions for conflict sets. 
\subsection{Geometrical description of different cases}
\begin{figure}[htbp]
 \centering   
\includegraphics*[width=0.4\textwidth]{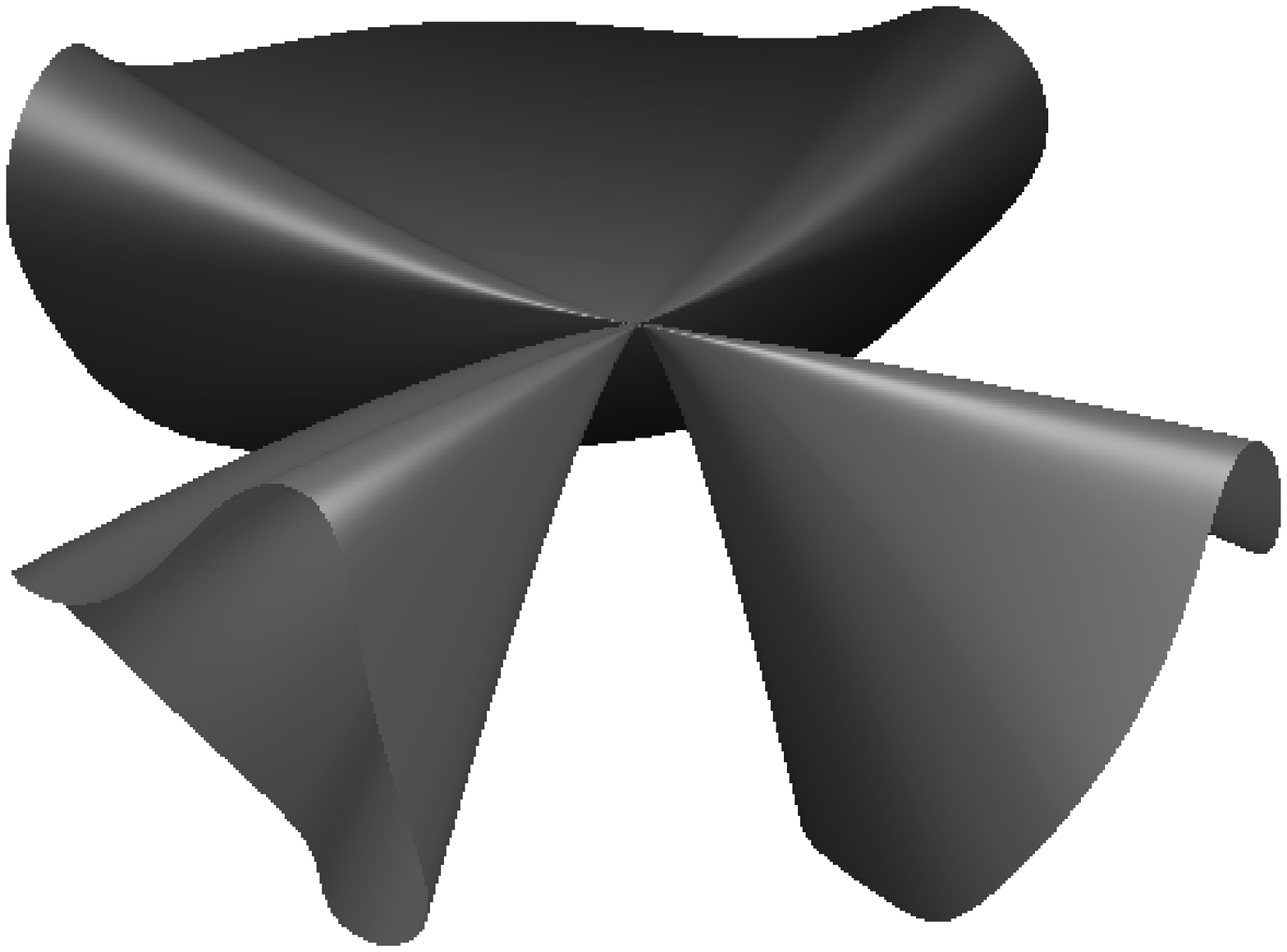}
\includegraphics*[width=0.4\textwidth]{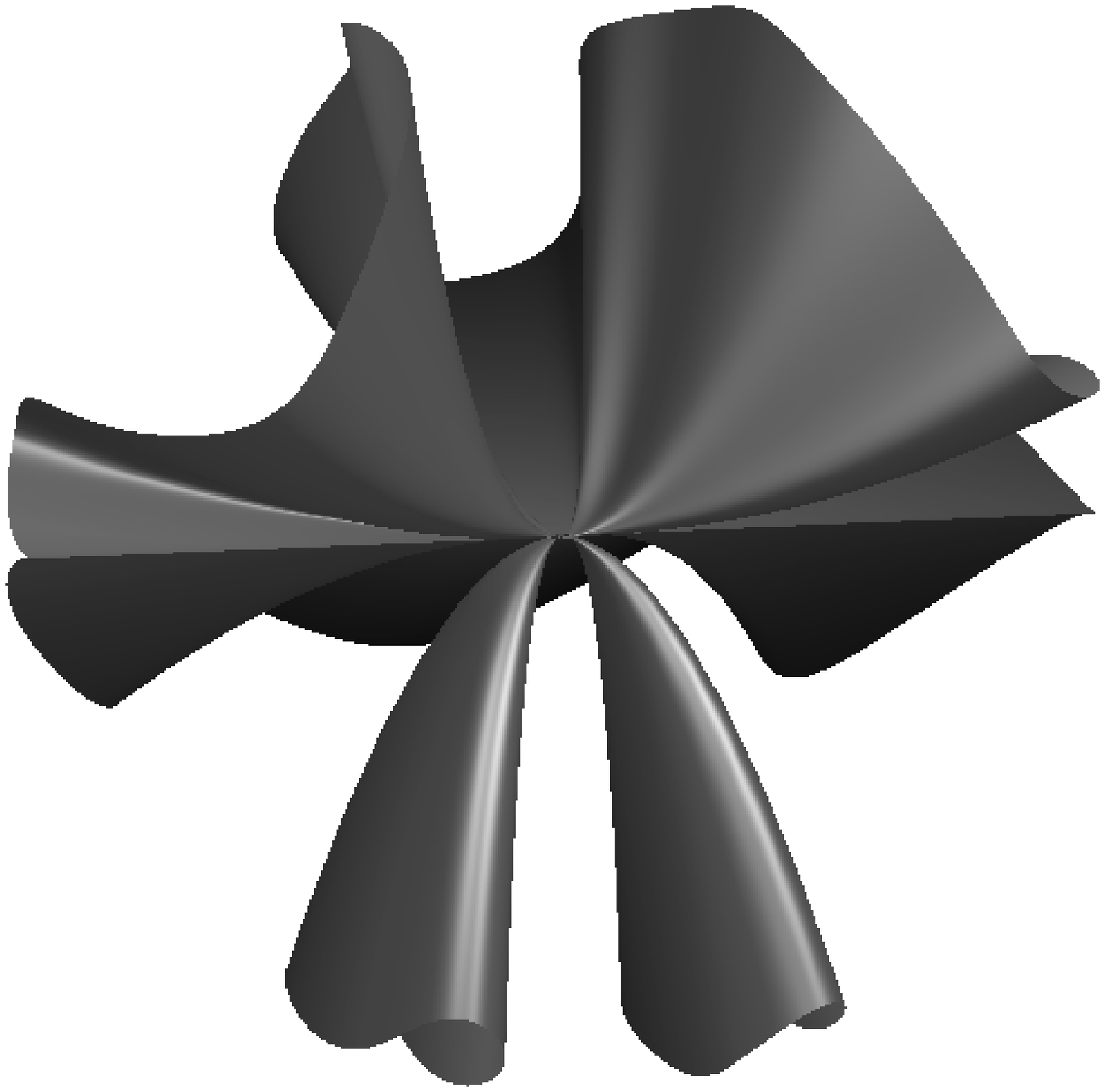}
\newline
\includegraphics*[width=0.4\textwidth]{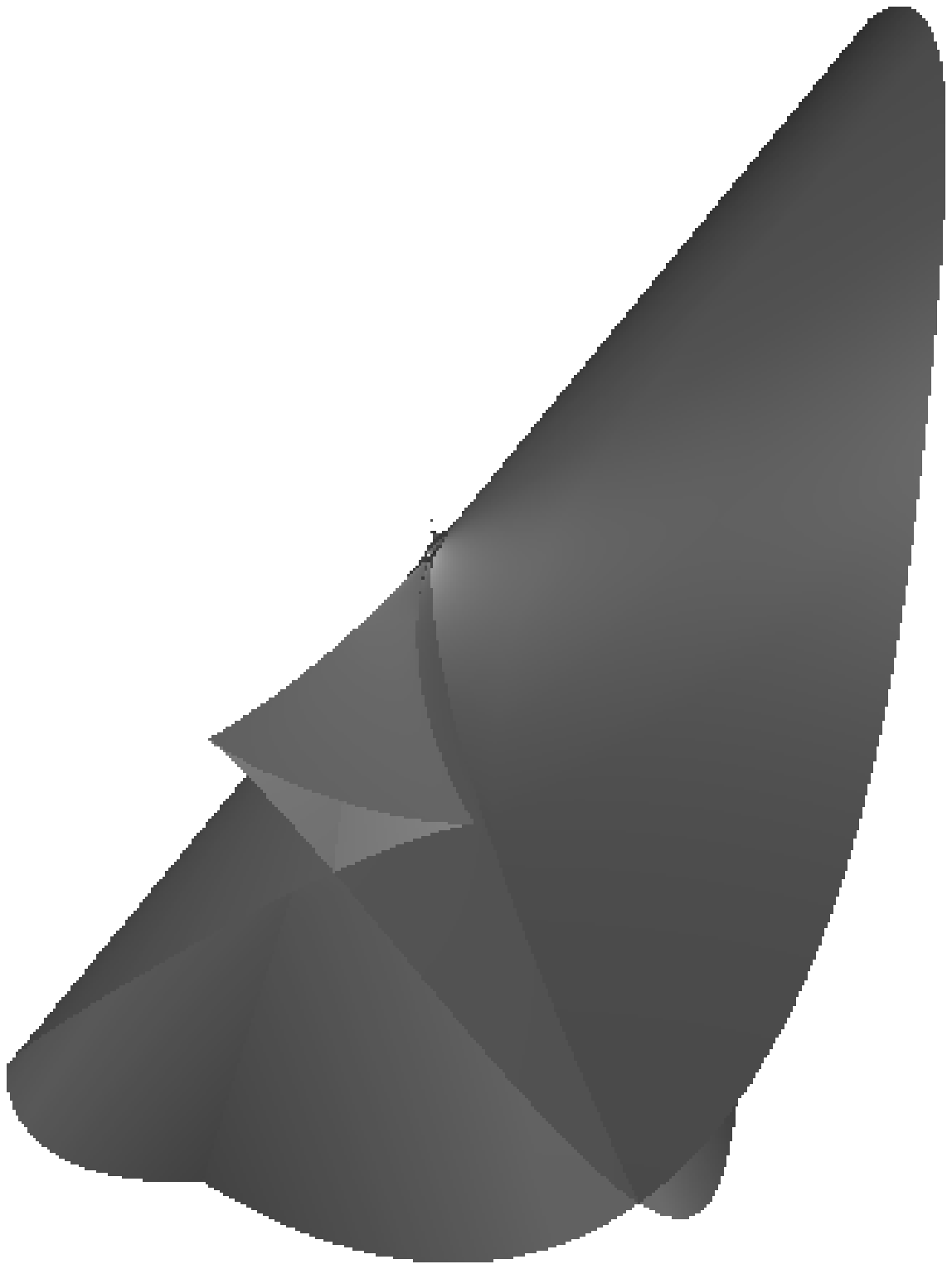}
\includegraphics*[width=0.4\textwidth]{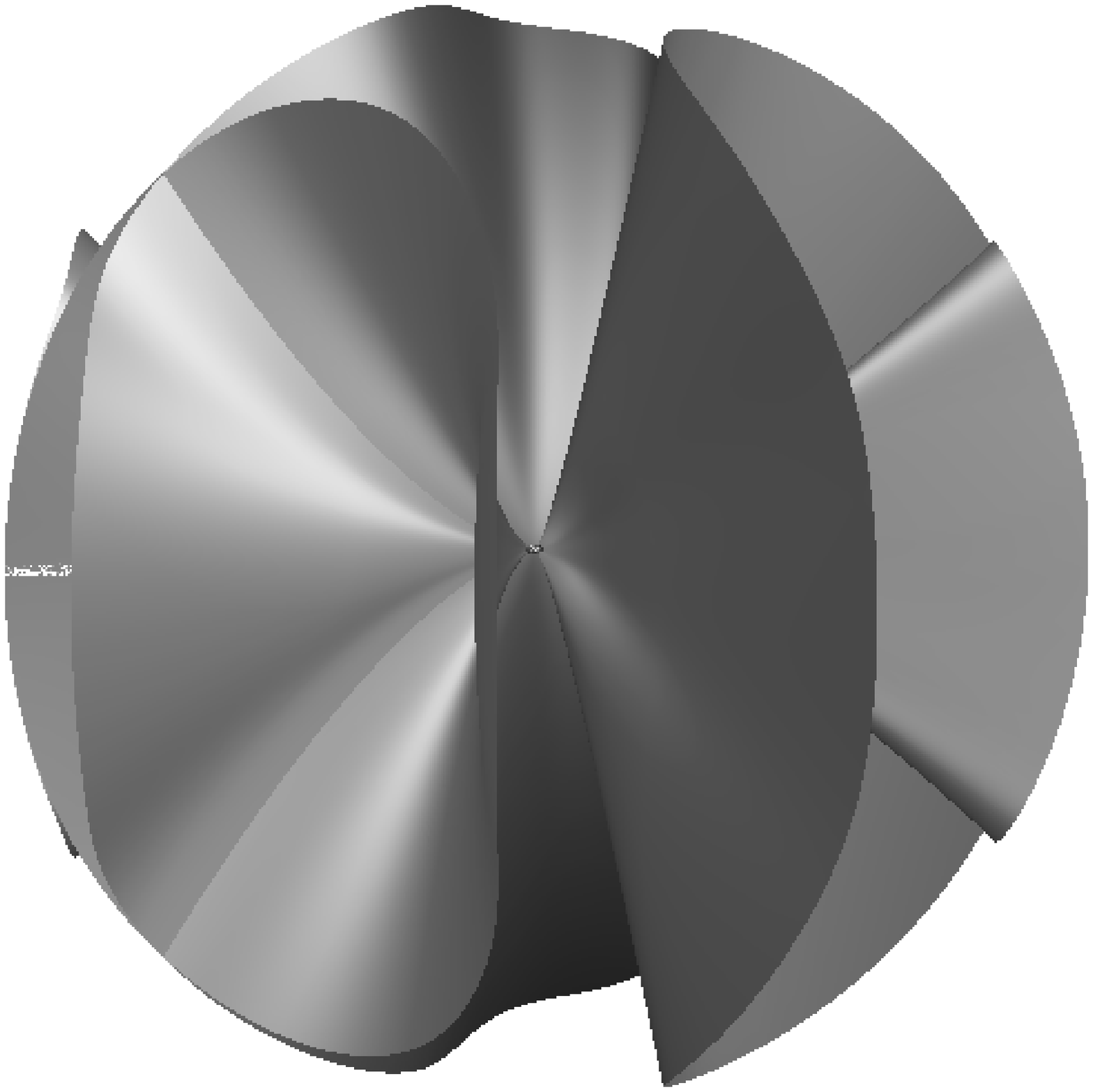}
    \caption{Sections of $A_3A_2$}
    \label{fig:a3a2}
\end{figure}

Once we know that the stratified big wavefronts intersect
transversally to determine what sort of singularities
can occur will follow from a codimension count. 
\newline
For the description of these singularities the main distinction is
the difference $n-l$.
\newline
Indeed if
$(\mu_1,\mu_2 , \cdots , \mu_l)$ is the list of codimensions
then we seek $\mu_i$ with $1 \leq \mu_i$ and $\sum_{i=1}^l\mu_i \leq n+1$.
Those $\mu_i$ that are $1$ correspond to smooth hypersurfaces. They are not
very interesting because they present just a reduction of
$n$ and $l$ by $1$, because if say $\mu_l=1$ then $N^*M_l^h$ is
smooth, and locally smoothly equivalent to $\R^{n-1}\times\R$.
Thus the singularity type reduces to what happens in 
$N^*M_1^h$, and thus it reduces to a problem with $l-1$ surfaces
in $\R^{n-1}$.
\newline
If $n-l$ is fixed then for arbitrary $n$ a certain number of
parts in the partition have to be $1$. Let $k$ be the number of strata that 
have codimension $>1$. It follows that $2k+(l-k) \leq n+1$ so that 
a maximum of $n-l+1$ codimensions is $>1$. The others are $1$.
\begin{prettylist}
\item[\fbox{$n-l=0$}]
If $n=l$ then at most $1$ of the $\mu_i$ is $>1$. So the only case
to consider is $l=2$. We can have only two cases: $(1)$, $(2)$.
\item[\fbox{$n-l=1$}]
At most $2$ of the codimensions are $>1$. So it is
enough to consider $n=3$, $l=2$. In addition to the above combinations
we will have: $(2,2)$ and $(3)$.
\item[\fbox{$n-l=2$}]
The relevant dimensions are: $n=5$, $l=3$. The new cases
are: $(4)$, $(3,2)$ and $(2,2,2)$.
\item[\fbox{$n-l=3$}]
Dimensions: $n=7$, $l=4$. 
New cases: $(5)$, $(4,2)$, $(3,3)$, $(3,2,2)$, $(2,2,2,2)$
\item[\fbox{$n-l=4$}]
Dimensions: $n=9$, $l=5$. 
New cases: $(6)$, $(5,2)$, $(4,3)$, $(4,2,2)$, $(3,3,2)$, $(3,2,2,2)$ and
$(2,2,2,2,2)$.
\end{prettylist}
For each of the strata there are only a limited number of singularities,
from the ADE list. The conflict set has dimension $n-l+1$. The codimension
of a singularity on a generic front of dimension $n-l+1$ is maximally
$n-l+2$. If we look at the above list we see that on the conflict set
the codimension can add up to $2(n-l+1)$. Thus the singularities we
encounter are the ones that we also expect to find in $n-l$ parameter
families of $n-l+1$ dimensional fronts, though this last list will
typically contain more singularities. 
\newline
One might ask whether the above multi-singularities do not
present any moduli. This is not the case.
Even though our singularities are not, if $n-l>0$, singularities of
generic fronts we are still allowed to produce local models - as
is remarked in \cite{BruceGiblinGibson} for the case $n=3$, $l=2$ -
because they present $R^+$-versal unfoldings of multigerms. 
\newline
So if $n-l=0$ the singularities of the conflict set are the
generic singularities of 2-dimensional fronts.
\newline
If $n-l=1$ the codimension can add up to $4$. The cases to consider
are $A_2A_2$, $A_1^2A_2$, $A_1^2A_1^2$. All other singularities are 
just those of generic $2$-dimensional fronts. Pictures are
again in \cite{BruceGiblinGibson}, but we take some time to
discuss a nice example. 
\newline
The $A_2A_2$ singularity is a generic projection of two
transversely intersecting cuspidal edges in $\R^4$. To obtain a picture
of this we take two copies of our previous example 
\begin{equation*}
  G_1\colon x_0=s_1^3+x_1s_1+x_2 \qquad G_2\colon x_0=s_2^3+x_3s_2-x_2
\end{equation*}
At zero these two intersect transversally. Next we 
project the intersection along the time axis $x_0$ to $\R^3$. 
The surface we get is the following picture.
\begin{figure}[htbp]
\centering
\includegraphics*[width=6cm]{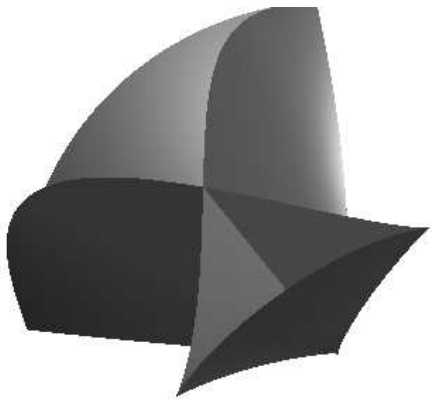}
\caption{The $A_2A_2$ surface}
\end{figure}
This is also known as $D_4^+$ if we view it as a metamorphosis
of a wavefront in $\R^3$. Recall that a metamorphosis is a 
one dimensional family of fronts, see\cite{MR93b:58019}.
The name $D_4^+$ is chosen because the surface is 
also obtained with an unfolding
\begin{equation*}
G_1-G_2=s_1^3-s_2^3+x_1s_1-x_3s_2+2x_2
\end{equation*}
This is not a versal unfolding. If we want to unfold the $D_4^+$ germ
$s_1^3-s_2^3$ with V-versal unfolding we need 4 parameters.  
\newline
If $n-l=2$ we need at least $n=4$ and $l=2$ to 
obtain an interesting new local model. Indeed the case $(4)$ has 
$A_4$ and $D_4^\pm$ and suspensions of the cases that occur
with $n-l=1$. So the first really new case is $(3,2)$. On this
stratum we have amongst others $A_3A_2$. This is a metamorphosis of
a 3-dimensional front. Some sections of this surface are in figure
\ref{fig:a3a2}. In one them we see a swallowtail meeting a cuspidal
edge.
\begin{rem}
All pictures here were obtained with the help of the software 
\cite{GPS01} and the program ``\texttt{surf}'', written by Stephan Endrass.
\end{rem}
\section{Concluding remark}
Similar results should hold for the center set. We expect that
if $n\leq 6$ the center set only has the multi-singularities
of wavefronts in that dimension. We also expect that smooth boundaries of
strictly convex compact domains in $\R^n$ the center symmetry
set should generically have only ADE-sings if $n\leq 6$. That is
because on such a boundary points with parallel tangent planes
always stay at a distance from each other. This distance is
uniformly bounded for such a surface, because of the
compactness. So they stay away from the diagonal and the local
situation can be treated as if the patches round the two
points with parallel tangent planes originated from two
distinct surfaces. Again this is all conjectural, though 
Giblin and Holtom have proved some assertions in this direction, see
\cite{MR2001e:58054}. 
\newcommand{\cprime}{{$'$}}
\providecommand{\bysame}{\leavevmode\hbox to3em{\hrulefill}\thinspace}


\begin{thebibliography}{AGZV85}

\bibitem[AGZV85]{Arnold:SAD}
V.~I. Arnol{\cprime}d, S.~M. Guse{\u\i}n-Zade, and A.~N. Varchenko,
  \emph{Singularities of differentiable maps. {V}ol. {I}}, Birkh\"auser Boston
  Inc., Boston, MA, 1985, The classification of critical points, caustics and
  wave fronts, Translated from the Russian by Ian Porteous and Mark Reynolds.

\bibitem[Arn90]{MR93b:58019}
V.~I. Arnol{\cprime}d, \emph{Singularities of caustics and wave fronts}, Kluwer
  Academic Publishers Group, Dordrecht, 1990.

\bibitem[BGM82]{MR84a:58001}
Thomas Banchoff, Terence Gaffney, and Clint McCrory, \emph{Cusps of {G}auss
  mappings}, Pitman (Advanced Publishing Program), Boston, Mass., 1982.

\bibitem[Dui74]{Duistermaat:CPAM}
J.\~J.\ Duistermaat, \emph{Oscillatory integrals, lagrange immersions and
  unfolding of singularities}, Comm. Pure Appl. Math. \textbf{27} (1974),
  207--281.

\bibitem[GH99]{MR2001e:58054}
Peter Giblin and Paul Holtom, \emph{The centre symmetry set}, Geometry and
  topology of caustics---CAUSTICS '98 (Warsaw), Polish Acad. Sci., Warsaw,
  1999, pp.~91--105.

\bibitem[GPS01]{GPS01}
G.-M. Greuel, G.~Pfister, and H.~Sch\"onemann, \emph{{\sc Singular} 2.0}, {A
  Computer Algebra System for Polynomial Computations}, Centre for Computer
  Algebra, University of Kaiserslautern, 2001, {\tt
  http://www.singular.uni-kl.de}.

\bibitem[H{\"o}r85]{Hormander:ALDO}
L.~H{\"o}rmander, \emph{The analysis of linear differential operators, i-iv},
  Springer Verlag, Berlin, 1983-1985.

\bibitem[JB85]{BruceGiblinGibson}
C.G.~Gibson J.W.~Bruce, P.J.~Giblin, \emph{Symmetry sets}, Proc. Cam.
  \textbf{145} (1985), 207--281.

\bibitem[Por94]{Porteous:1994}
I.~Porteous, \emph{Geometric differentiation}, C.U.P., 1994.

\bibitem[Sie99]{Siersma1}
Dirk Siersma, \emph{Properties of conflict sets in the plane}, Geometry and
  topology of caustics---CAUSTICS '98 (Warsaw), Polish Acad. Sci., Warsaw,
  1999, pp.~267--276.

\bibitem[Wal77]{Wall}
C.T.C.\ Wall, \emph{Geometric properties of generic differentiable manifolds},
  Geometry and topology (Proc. III Latin Amer. School of Math., Inst. Mat. Pura
  Aplicada CNPq, Rio de Janeiro, 1976), Springer, Berlin, 1977, pp.~707--774.
  Lecture Notes in Math., Vol. 597.

\end{thebibliography}
\end{document}